\documentclass[12pt]{article}

\def\Endproof{{\ \vbox{\hrule\hbox{%
   \vrule height1.3ex\hskip0.8ex\vrule}\hrule
}}}
\def\0{\leqno}
\author{{\small by}\\Gabriel Mitric and Izu Vaisman}
\date{}
\title{POISSON STRUCTURES ON TANGENT BUNDLES
}
\begin{document}
\maketitle
{\def\thefootnote{*}\footnotetext[1]%
{{\it 2000 Mathematics Subject Classification}: 53D17.
Graded bivector field.}}
\begin{center} \begin{minipage}{12cm} A{\footnotesize BSTRACT.
The paper starts with an interpretation of the complete lift of a
Poisson structure from a manifold $M$ to its tangent bundle $TM$
by means of the Schouten-Nijenhuis bracket of covariant symmetric
tensor fields defined by the cotangent Lie algebroid of $M$. Then,
we discuss Poisson structures of $TM$ which have a graded
restriction to the fiberwise polynomial algebra; they must be
$\pi$-related ($\pi:TM\rightarrow M$) with a Poisson structure on
$M$. Furthermore, we define {\it transversal Poisson structures}
of a foliation, and discuss bivector fields of $TM$ which produce
graded brackets on the fiberwise polynomial algebra, and are
transversal Poisson structures of the foliation by fibers.
Finally, such bivector fields are produced by a process of {\it
horizontal lifting} of Poisson structures from $M$ to $TM$ via
connections.}
\end{minipage}
\end{center} \noindent
\section*{\begin{center}
{\bf 1. The complete lift of a Poisson structure}
\end{center}}
Let $M$ be an $n-$dimensional differentiable manifold with  local
coordinates $(x^i),$  $(i=1,...,n),$  $\pi :TM\longrightarrow M$
its tangent bundle and $(y^i)$ the vector coordinates with respect
to the basis $\{ {\partial}/{\partial x^i}\}.$ (We assume that
everything is $C^{\infty}$ in this paper).

Let us consider a Poisson structure on the manifold $M,$ given by
the Poisson bivector
$$
w=\frac {1}{2}w^{ij}\frac {\partial}{\partial x^i}\wedge \frac
{\partial}{\partial x^j}
 \0(1.1)
$$
(throughout the paper, we use Einstein's summation convention).
The complete lift of $w$ in the sense of \cite{KY} is given by
$$
w^C=w^{ij}\frac {\partial}{\partial x^i}\wedge \frac
{\partial}{\partial y^j}+ \frac {1}{2}y^k \frac {\partial w^{ij}
}{\partial x^k}\frac {\partial}{\partial y^i}\wedge \frac
{\partial}{\partial y^j} \ ,
 \0(1.2)
$$
and it follows easily that $w^C$ is a Poisson bivector field on
$TM$ since the Poisson condition, namely, that the
Schouten-Nijenhuis bracket $[w^C,w^C]=0$ \cite{Va1}, is satisfied.

The Poisson bivector $w^C$ has already been studied by several
authors  \cite{{Va3},{Wein1}}, and  it can also be  derived
 from the bracket of the $1-$forms of $M$ with respect to
the Poisson structure $w$ (e.g., \cite{Va1})
$$
\{\alpha, \beta \}= {\mathcal{L}}_{{\alpha }^{\sharp}}\beta
-{\mathcal{L}}_{{\beta }^{\sharp}}\alpha -d(w(\alpha,\beta)) \ \ \
(\alpha ,\beta \in \Omega ^1(M)).
 \0(1.3)
$$

A Pfaff form $\alpha = \alpha _idx^i$ on $M$ may be regarded as a
fiberwise  linear function $l(\alpha ):=\alpha _i(x)y^i$ on $TM$
($:=$ denotes a definition).
 A Poisson structure $W$ on $TM$ is completely determined by the
brackets $\{f\circ \pi, g \circ \pi\}_W,$ $\{l(\alpha ), f\circ
\pi \}_W$ and $\{l(\alpha ),l(\beta )\}_W,$ where  $ f,g \in
C^{\infty}(M)$ and  $\alpha ,\beta \in {\Omega}^1(M),$ the space
of Pfaff forms on $M$, since it is completely  determined by
 the brackets of the local coordinates $x^i$ and  $y^j$.

The Poisson bivector $w^C$ is exactly the one defined by:

$i) \ \{f\circ \pi, g \circ \pi\}_{w^C} =0 \ , \ \ \ \forall f,g
\in C^{\infty}(M);$

$ii) \ \{l(\alpha ), f \circ \pi\}_{w^C} =({\alpha
}^{\sharp}f)\circ \pi, \ \ \forall f \in C^{\infty}(M), \  \forall
\alpha \in \Omega ^1(M),$
\newline
where  $\sharp _{w}:T^{*}M\longrightarrow TM$ is defined by
 $ \beta({\alpha}^{\sharp } )=w(\alpha , \beta ), \ $ $\forall \beta
\in {\Omega}^1(M); $

$iii)  \ \{l(\alpha ),l(\beta )\}_{w^C} =l(\{\alpha,\beta\}) \ , \
\ \forall \alpha , \beta  \in {\Omega }^1(M)$.

$w^C$ is a Poisson structure  because the bracket (1.3) satisfies
the Jacobi identity.

The Poisson structure $w^C$ also has the interesting property
$$w^C=-{\mathcal{L}}_Ew^C, \hspace{5mm} E=y^i\frac{\partial}{\partial
y^i}$$ ($E$ is the {\it Euler vector field}), which means that
$(TM,w^C)$ is a {\it homogeneous Poisson manifold} \cite{Va3}.

\vspace{0.3cm} \hspace{0.3cm} We remind that a {\it Lie algebroid
} over a manifold $M$ is a triple $(A,{[ \ , \ ]}_A,\sigma ),$
where $p:A\longrightarrow M$ is a vector bundle, ${[ \ , \ ]}_A$
is an $\mathbf{R}-$Lie algebra structure on the space $\Gamma A$
of the global cross sections of $A$
 and $\sigma :A\longrightarrow TM$ is a  morphism of
vector bundles, called {\it anchor,} such that
$$
i) \ \sigma ([s_1,s_2]_A)=[\sigma (s_1), \sigma (s_2)],\; ii) \
[s_1,fs_2]_A=f[s_1,s_2]_A+((\sigma s_1)f)s_2
$$
for every $s_1,s_2 \in \Gamma A$, $f \in C^{\infty}(M)$, and where
$[\;,\;]$ is the Lie bracket of vector fields on $M$.

In what follows, we give one more interpretation of the Poisson
structure $w^C$ by means of a Schouten-Nijenhuis bracket on a Lie
algebroid $A.$

There exists a well known operation, called the {\it
Schouten-Nijenhuis bracket}, on cross sections of
 ${\mathcal V}(A):={\oplus }_k\Gamma{\wedge }^kA$ (e.g., see \cite{YKS1}).
 A less popular
 operation, the {\it Schouten-Nijenhuis bracket of symmetric
 tensors} also exists, and was studied in an algebraic context and
 for $TM$ \cite{Ba}. Here we present this second operation on
 the algebra of cross sections
 $S(A)=\bigoplus \limits _{k\geq 0}S_k(A)$, where
$S_k(A)= \Gamma {\odot}^kA,$ $A$ is a Lie algebroid,  $\odot $
denotes the symmetric tensor product, and $\Gamma$ denotes spaces
of global cross sections of bundles. Then, we show that this
operation leads to another definition of the complete lift $w^C.$

For any Lie algebroid $A,$ one has the {\it Lie derivative}
\cite{Ba} which is defined
 by putting ${\mathcal{L}}^A_s f={\mathcal{L}}_{\sigma (s)}f$
for functions $f\in C^{\infty}(M),$ and
${\mathcal{L}}^A_ss^{'}={[s,s^{'}]}_A$ for cross sections
$s^{'}\in \Gamma A,$ and by  extending it to arbitrary cross
sections of $\Gamma (({\otimes}^kA)\otimes ({\otimes}^lA^{*}))$ as
in the case of the classical Lie derivative.
 In particular, we have the restriction $
{\mathcal{L}}_{s}^A: S_{k}(A) \longrightarrow S_{k}(A)  .
 $

 PROPOSITION 1.1.  {\it There exists a well defined unique
  extension of local type  of the Lie derivative ${\mathcal{L}}^A_{s}$ to an
$\mathbf{R}-$bilinear operation
$$
< \ , \ > :S_{p}(A)\times S_{q}(A) \longrightarrow S_{p+q-1}(A)  \
\ \ (p,q\geq 1),
$$
such that
$$
<s_1\odot ... \odot s_p, t_1 \odot ... \odot t_q>=
\0(1.4)
$$
$$=\sum _{i=1}^p\sum _{j=1}^q[s_i,t_j]\odot s_1\odot
...\odot\hat{s_i}\odot ... \odot s_p \odot t_1 \odot ... \odot
\hat{t_j} \odot ... \odot t_q \ ,
$$
where the hat denotes the absence of the corresponding factor.}

{\it Proof.}  ``Local type"  means that  $\forall x_0\in M$ and
$\forall G\in S_p(A),$ $H\in S_q(A),$ $<G,H>(x_0)$ depends only on
the restriction of the cross sections  $G,$ $H$ to a neighborhood
of $x_0.$

Since, around $x_0\in M$ $G,H$ are decomposable into finite sums
of products of the form appearing in $(1.4),$ the required
extension has an obvious definition, and we only must show that
the result does not depend on the decomposition.

 If $H=t_1 \odot ... \odot t_q $, formula $(1.4)$ becomes
$$
<s_1\odot ... \odot s_p, H>=\sum
_{i=1}^p({\mathcal{L}}^A_{s_i}H)\odot s_1\odot
...\odot\hat{s_i}\odot ... \odot s_p  \ .
 \0(1.5)
$$
Similarly, if $G=s_1 \odot ... \odot s_p$ then
$$
<G, t_1\odot ... \odot t_q>=-\sum
_{j=1}^q({\mathcal{L}}^A_{t_j}G)\odot t_1\odot
...\odot\hat{t_j}\odot ... \odot t_q  \ .
 \0(1.6)
$$
Formulas  $(1.5),$  $(1.6)$ show the required independence of the
bracket $<G,H>$ of the decomposition of $G,$ $H$.  \Endproof

Note that we may also  consider $(1.5),$  $(1.6)$ as the
definition of the bracket $<G,H>$.

Now, we may extend the bracket $< \ , \ >$ to the case where the
factors belong to $S_0(A)=C^{\infty}(M).$ Namely, we will put
$$
<g,h>=0, \ \ <s_1\odot ... \odot s_p, f>= - <f,s_1\odot ... \odot
s_p>  \0(1.7)$$
$$=
 \sum _{i=1}^p({\mathcal{L}}^A_{s_i}f) s_1\odot
...\odot\hat{s_i}\odot ... \odot s_p \ ,
$$
$\forall f,g,h\in C^{\infty }(M),$ $s_i\in \Gamma A.$ The fact
that the second formula $(1.7)$ does not depend on the
decomposition $G=s_1\odot ... \odot s_p$ follows by noticing that
$$
{\mathcal{L}}^A_{s_i}f =(d_Af)(s_i) \ ,
$$
where $d_A$ is the exterior differential for the Lie algebroid $A$
 \cite{YKS1}, whence
$$
<G,f>=i(d_Af)G,
 \0(1.8)
$$
where the definition used for the operator $i$ is that of
\cite{KN}.

 PROPOSITION 1.2. {\it The bracket  $<,>$   has the
following properties}
$$
<H,G>=-<G,H> \ ,
 \0(1.9)
$$
$$
<G,H\odot K>=<G,H>\odot K+H\odot <G,K> \ ,
 \0(1.10)
$$
 $\forall \  G,H,K \in S(A).$

{\it Proof.} The bracket $<,>$ is extended to $S(A)$ by
$\mathbf{R}$-bilinearity.
 Both relations  easily follow from $(1.4)$ and $(1.7).$
  \Endproof

 PROPOSITION 1.3. (The Jacobi identity.)  $\forall F,G,H \in
 S(A)$,
$$
<<G,H>,K>+<<H,K>,G>+<<K,G>,H> \ =0 \ .
\0(1.11)
$$

{\it Proof.} It suffices to prove $(1.11)$ for decomposable $G,$
$H,$ $K$ and this follows by a technical computation bases on the
Jacobi identity satisfied by the bracket of the cross sections of
$A.$ \Endproof

 COROLLARY 1.4. {\it $(S(A),<  ,  >)$ is a Poisson algebra
\cite{Ba} with respect to the symmetric product $\odot $ and the
Schouten-Nijenhuis bracket}  $<  ,  >.$

\vspace{2mm}
 Let us consider  the particular case
 of the
cotangent Lie algebroid $(T^{*}M, \{ \ , \ \},$ $ {\sharp}_w)$ of
a Poisson manifold $(M,w),$ where the bracket is that defined by
formula (1.3).

 Then
$$
S(T^{*}M):=\oplus _k S_k(T^{*}M)\leqno{(1.12)}
$$
is the algebra of the covariant symmetric tensor fields on $M,$
and Corollary 1.4 shows that
  $(S(T^{*}M),<  ,  >)$ is a Poisson algebra with respect to
the symmetric product $\odot$  and the Schouten-Nijenhuis bracket
$<  ,  >$ of the Lie algebroid $T^{*}M.$

Notice that, in the present case $d_{T^{*}M}f=-X_f,$ where $X_f$
is the $w-$Hamil\-tonian vector field of $f.$ Accordingly, $(1.8)$
yields
$$
<G,f>=-i_{X_f}G \ \ \  \ ( f\in C^{\infty}(M), \  G\in
S_p(T^{*}M)).
 \0(1.13)
$$

The function space $C^{\infty}(TM)$ has some interesting
subspaces. Namely, the spaces of  fiberwise homogeneous
$k-$polynomials
$$
{\mathcal{H}\mathcal{P}}_k(TM):=\{\tilde{G}=G_{i_1...i_k}y^{i_1}...y^{i_k}
/
  \0(1.14)
$$
$$
  G=G_{i_1...i_k}dx^{i_1}\odot ... \odot dx^{i_k}\in
 S_k(T^{*}M)\} \ ,
$$
and we have an isomorphism of algebras
 $$
\iota \ : (S(T^{*}M),\odot )\longrightarrow {(\mathcal{H}
\mathcal{P}}(TM):={\oplus }_k{\mathcal{H}}{\mathcal{P}}_k(TM), \
\cdot \ )
  \0(1.15)
 $$
mapping $G$ to $\tilde{G}:=\iota(G)$ (the dot denotes usual
multiplication). With this isomorphism, the bracket $< , >$ of
symmetric covariant tensor fields is translated into a bracket of
polynomials. Moreover, since the local coordinates $x^i,$ $y^j$
are polynomials of degree zero and one respectively, this bracket
defines a Poisson structure with a Poisson bivector, say $W$, on
$TM$; the bracket will be denoted by $\{\;,\;\}_W$.

 PROPOSITION 1.5. {\it The Poisson structure $W$ defined on the
 tangent bundle of a Poisson manifold  $(M,w)$ by the bracket
 $\{\;,\;\}_W$
  coincides with the Poisson structure $w^C.$ }

{\it Proof.} The  brackets $\{x^i,x^j\}_W,$ $\{x^i,y^j\}_W$ and
$\{y^i,y^j\}_W$ computed with  $< , >,$ are the same as those
produced by $(1.2).$ \Endproof

 COROLLARY 1.6. {\it   If $G$ and $H$ are symmetric covariant
tensor fields on $M,$ then
$$
\widetilde{<G,H>}=\{\tilde{G}, \tilde{H}\}_{w^C} \ ,
   \0(1.16)
$$
and $(1.15)$ is a corresponding isomorphism of Poisson algebras.}
\vspace{0.3cm}

 In the remaining part of this section we will
compute the {\it modular class} of the Poisson structure $w^C$.

If $\mu $ is a volume form on an orientable manifold $M,$ the
divergence ${div}_{\mu}X$ of a vector field $X$ is defined by the
condition
$$
{\mathcal{L}}_X\mu=({div}_{\mu}X)\mu  \ ,
$$
and one has
$$
{div}_{\mu}(fX)=f{div}_{\mu}X \ + \ Xf \ , \ \ f\in C^{\infty}(M)
\ .
$$
Accordingly, if $(M,w)$ is a Poisson manifold endowed with a
volume form $\mu ,$  the operator
$$
{\Delta }_{\mu}:f\in C^{\infty}(M)\longmapsto
 {div}_{\mu}X_{f}\in
C^{\infty}(M)
$$
is a derivation on $C^{\infty}(M),$ so it is a vector field on the
manifold $M,$  called the {\it modular vector field } of
$(M,w,\mu)$ (see \cite{{YKS2},{W2}}).

Denote by ${\mathcal{V}}^i(M)$   the space of {\it $i$-vector
fields}
 of a manifold $M,$ i.e., skew symmetric contravariant tensor fields of
type $(i,0)$ on $M,$ and ${\mathcal{V}}(M)=(\bigoplus \limits
_{i=1}^n {\mathcal{V}}^{i}(M), \wedge )$  the {\it contravariant
Grassmann algebra} of $M.$ On a Poisson manifold $(M,w),$ the
Lichnerowicz-Poisson coboundary operator  is
$$
\sigma : =-[w, \ . \ ]:{\mathcal{V}}^{k}(M)\longrightarrow
{\mathcal{V}}^{k+1}(M) \ ,
$$
 where $[ \  , \ ]$ is the Schouten-Nijenhuis bracket,
 and one has the Lichnerowicz-Poisson (LP) cohomology spaces
(e.g., \cite{Va1})
$$
H^k_{LP}(M,w)=\frac {Ker \ (\sigma:{\mathcal{V}}
^k(M)\longrightarrow {\mathcal{V}}^{k+1}(M) )}{Im \
(\sigma:{\mathcal{V}} ^{k-1}(M)\longrightarrow
{\mathcal{V}}^{k}(M) )} \ .
$$

For a modular vector field one has $\sigma {\Delta} _{\mu}=0$
\cite{YKS2}, and ${\Delta} _{\mu}$ is a $1-$cocycle. Therefore it
defines a $1-$dimensional  LP-class $ \Delta =[\Delta _{\mu}] \in
H^1_{LP}(M,w).$ It is easy to see  that this class does not depend
on $\mu$; it is called the {\it modular class} of the Poisson
manifold $(M,w)$ \cite{{YKS2},{W2}}.

We want to discuss the relation between the modular classes of
$(M,w)$ and of $(TM,w^C).$

Let $g$ be a Riemannian metric on  the oriented manifold $M.$ Then
$$
dV_g=\sqrt{\det g} \ dx^1\wedge... \wedge dx^n
\0(1.17)
$$
 is a volume form on $M,$ and it follows easily that
$$
\Phi =(\det g) dx^1\wedge... \wedge dx^n\wedge dy^1 \wedge ...
\wedge dy^n
\0(1.18)
$$
is a volume form on $TM$ (the volume form of the Sasaki metric
associated to $g$ \cite{Dom}).

 PROPOSITION 1.7.  {\it The modular vector field of $(TM,w^C,\Phi
)$ is given by
$$
{\Delta}_{\Phi }^{TM}=2({\Delta}_{dV_g}^{M})^V,
\0(1.19)
$$
where the upper index $V$ denotes the vertical lift in the sense
of} \cite{KY} .

{\it Proof. } With $(1.1),$ a Hamiltonian vector field on $(M,w)$
has the form
$$
X_f^w=\{f, \cdot \}_{w}= \frac {\partial f}{\partial
x^i}w^{ij}\frac {\partial}{\partial x^j} \  \ \ \ (f\in
C^{\infty}(M)) \ ,
 \0(1.20)
$$
and the definition of the modular vector field leads to
$$
{\Delta}_{dV_g}=\sum \limits _{k=1}^{n} \left( \frac {\partial
w^{ik}}{\partial x^k} +w^{ik}\frac {\partial \ln \sqrt{\det
g}}{\partial x^k } \right)\frac {\partial }{\partial x^i} \ .
 \0(1.21)
$$

Then,  if $F\in C^{\infty}(TM),$
 $(1.2)$ gives for the
Hamiltonian vector field $X_F^{w^C}$ the expression
$$
X_F^{w^C}=\frac {\partial F}{\partial y^k}w^{ki}\frac
{\partial}{\partial x^i}+ \left( \frac {\partial F}{\partial
x^k}w^{ki}+\frac {\partial F}{\partial y^k}y^h \frac {\partial
w^{ki} }{\partial x^h}\right) \frac {\partial}{\partial y^i}
  \ ,
 \0(1.22)
$$
and a straightforward computation yields the modular vector field
$$
{\Delta}_{\Phi }=\sum \limits _{k=1}^{n} 2\left( \frac {\partial
w^{ik}}{\partial x^k} +  w^{ik}\frac {\partial \ln \sqrt{\det
g}}{\partial x^k } \right)\frac {\partial }{\partial y^i} \ .
 \0(1.23)
$$
This exactly is the required result. \Endproof

 COROLLARY 1.8. {\it The modular class of the Poisson manifold
$(TM,w^C)$ is represented by $2{\Delta}_{\mu}^V,$ for every
modular vector field ${\Delta}_{\mu }$ of the base manifold
$(M,w).$ }

{\it Proof.} Proposition 1.7 shows that the result is true for the
field ${\Delta}_{dV_g}.$ Since  from $(1.20),$ $(1.22)$ one also
gets
$$
({\sigma }_wf)^V= {\sigma }_{w^C}(f\circ \pi ), \ \ f\in
C^{\infty}(M),  \ \pi :TM \longrightarrow M,
$$
the result is true for any other modular vector field of $(M,w).$
\Endproof

It is also worthwhile to notice the following result.

 PROPOSITION 1.9. {\it   The complete lift of multivector fields
induces a homomorphism of cohomology algebras }
$$
[Q]\in H^{*}_{LP}(M,w) \longmapsto [Q^C]\in H^{*}_{LP}(TM,w^C) \ .
 \0(1.24)
$$

{\it Proof.} The complete lift of multivector fields is the
natural extension of the complete lift of vector fields, and is
compatible with the Lie bracket \cite{KY}. Therefore, since the
Schouten-Nijenhuis bracket extends the Lie bracket e.g.,
\cite{Va1},
 if $Q_1,Q_2\in {\mathcal{V}}(M),$
$$
[Q_1,Q_2]^C=[Q_1^C,Q_2^C] \ ,
$$
where the bracket is the Schouten-Nijenhuis bracket. This implies
$$
(\sigma _wQ)^C=\sigma _{w^C}Q^C \ ,
  \0(1.25)
$$
and it follows that $(1.24)$ is a homomorphism. \Endproof \newpage
\section*
{\begin{center} {\bf 2. Graded Poisson structures on tangent
bundles}
\end{center}}
Recall that on $TM$ we have the spaces of fiberwise polynomial
functions ${\mathcal{H}}{\mathcal{P}}_{k}$ given by $(1.14)$.
Denote by
$$
{\mathcal{P}}_{k}(TM):= \bigoplus \limits
_{h=0}^k{\mathcal{H}\mathcal{P}}_{h} \ , \ \ {\mathcal{P}}(TM):=
 \bigcup \limits _{k\geq 0}{\mathcal{P}}_{k} \ ,
$$
the space of fiberwise non homogeneous polynomials of degree $\leq
k,$ and the {\it polynomial algebra.}

In particular, we have  the space
${\mathcal{A}}(TM):={\mathcal{P}}_{1}(TM)$ of  {\it affine
functions}
$$
a(x,y)=f(x)+l(\alpha ) \ , f\in C^{\infty}(M), \ \alpha \in \Omega
^1(M),
$$ where $l(\alpha)$ was defined in Section 1,
 and the space ${\mathcal{P}}_2(TM)$
 of non-homogeneous quadratic polynomials:
$$
p(x,y)=f(x)+l(\alpha )+s(G)
$$
where $G=G_{ij}dx^i\odot dx^j$ is a symmetric covariant tensor
field on $M$ and $s(G):=\tilde{G}$ is defined in $(1.14).$ Here
and in the whole paper, when speaking of polynomials on $TM,$ we
always mean fiberwise polynomials.

 DEFINITION 2.1. A Poisson structure $W$ on $TM$ is called {\it
polynomially graded} if ${\mathcal{P}}(TM)$ is closed by Poisson
brackets and $\forall F,G\in {\mathcal{P}}(TM)$
$$
F\in {\mathcal{P}}_{h}, \ G\in {\mathcal{P}}_k \Longrightarrow
\{F,G\}_W\in {\mathcal{P}}_{h+k} \ .
  \0(2.1)
$$

 PROPOSITION 2.2. {\it A polynomially graded Poisson structure
$W$ on $TM$ induces a Poisson structure $w$
 on the base manifold $M,$ such that the projection
 $\pi :(TM,W)\longrightarrow (M,w)$ is a Poisson mapping.   }

{\it Proof.} If  $f\in C^{\infty}(M),$  $f=f\circ \pi$ is a
polynomial  of degree zero on $TM.$ Thus, by $(2.1),$ $\forall f,g
\in C^{\infty}(M)$, and
$$
\{f,g\}_w:=\{f\circ \pi, g\circ \pi \}_W, \0(2.2)
$$
defines a Poisson structure $w$ on $M.$ \Endproof

Hereafter, we write $f$ for both  $f\in C^{\infty}(M),$ and
$f\circ \pi \in C^{\infty}(TM).$  The bracket $\{ \ , \ \}_W$ will
be denoted simply by $\{ \ , \ \}.$

Proposition 2.2 tells us that the polynomially graded Poisson
structures $W$ of $TM$ (if any) are {\it lifts} of Poisson
structures $w$ of $M$ i.e., $\pi:(TM,W)\rightarrow(M,w)$ is a
Poisson mapping. We suggest that the general problem of looking
for lifts of Poisson structures of a manifold to its tangent
bundle is an interesting problem.

The polynomially graded Poisson structure $W$ is completely
determined if, along with the brackets $\{f,g\},$ we also define
the brackets $\{l (\alpha),f\}$ and $\{l (\alpha),l (\beta )\},$
where $\alpha , \beta \in \Omega ^1(M).$

By $(2.1),$ the bracket $\{l (\alpha),f\}$ $\in $
${\mathcal{P}}_1(TM)$ i.e.,
$$
\{l(\alpha),f\}=X_{\alpha}f+l(\gamma _{\alpha}f) \ ,
\0(2.3)
$$
where $X_{\alpha}f\in C^{\infty}(M)$ and $\gamma _{\alpha}f \in
\Omega ^1(M).$

Since $\{l(\alpha),\ . \ \}|_{ C^{\infty}(M)}$ is a derivation of
$C^{\infty}(M), \ $
 it follows that
 $ X_{\alpha} $
is  a vector field on $M,$ and the mapping $
  \gamma _{\alpha }:C^{\infty}(M)\longrightarrow \Omega ^1(M)
$ also  is a derivation. Therefore, $\gamma _{\alpha}f $
 only depends  on $df.$

 The Leibniz rule implies
$$
\{l(h\alpha ) ,f\}=h(X_{\alpha }f)+l((X_h^{w}f)\alpha +h(\gamma
_{\alpha}f)) \ .
$$
Hence $\gamma$ must satisfy
$$
 \gamma _{h\alpha}f=h\gamma
_{\alpha}f+ (X_h^{w}f)\alpha \ .
\0(2.4)
$$

Similarly, the  bracket $\{l(\alpha),l(\beta )\}$
 must have an  expression of the form
$$
\{l(\alpha),l(\beta )\}=U(\alpha , \beta )+l(\Phi (\alpha ,\beta
))+s(\Psi (\alpha ,\beta )) \ ,
\0(2.5)
$$
where
$$
U(\alpha , \beta )\in C^{\infty}(M), \ \Phi (\alpha , \beta) \in
\Omega ^1(M), \;     \Psi (\alpha , \beta) \in S_2(T^{*}M)
$$
 are skew-symmetric operators.
 A replacement of
 $\beta $ by $f\beta  \ $  in $(2.5)$ leads to
 $$
U(\alpha ,f \beta )=fU(\alpha , \beta ) \ ,
 $$
whence,  $U$ is a bivector field on $M,$ and
$$
 \Phi (\alpha , f\beta )=f\Phi
(\alpha , \beta )+(X_{\alpha }f)\beta \ , \ \  \Psi (\alpha ,
f\beta )=f\Psi (\alpha , \beta )+\gamma_{\alpha }f \odot \beta \ .
\0(2.6)
$$

On the other hand, if $(x^i)$ are local coordinates on $M,$ if
$y^i=l(dx^i),$ and if $w$ is the Poisson structure introduced by
Proposition 2.2, Definition 2.1 tells us that the local coordinate
expression of $W$ must be of the form
$$
W=\frac {1}{2}w^{ij}(x)\frac {\partial}{\partial x^i}\wedge \frac
{\partial}{\partial x^j}+ (\varphi ^{ij}(x)+y^aA^{ij}_a(x))\frac
{\partial}{\partial x^i}\wedge \frac {\partial}{\partial y^j}
\0(2.7)
$$
$$
+\frac {1}{2}(\eta ^{ij}(x)+y^a\chi
^{ij}_a(x)+y^ay^bB^{ij}_{ab}(x))\frac {\partial}{\partial
y^i}\wedge  \frac {\partial}{\partial y^j} \ ,
$$
where $w,$ $\varphi ,$ $A,$ $\eta ,$ $\chi ,$ $B$ are local
functions on $M.$

 DEFINITION 2.3. A  polynomially graded Poisson structure $W$ on
$TM$ is said to be a {\it graded structure} if $\forall
 F\in {\mathcal{H}\mathcal{P}}_{h},$ $\forall
 G\in {\mathcal{H}\mathcal{P}}_{k},$
 $\{F,G\}_W\in {\mathcal{H}\mathcal{P}}_{h+k} \ .$

The conditions for a polynomially graded structure on $TM$ to be
graded are $ X_{\alpha }=0, \ $ $ U=0,$ $\Phi =0,$ and then
$(2.7)$ reduces to
$$
W=\frac {1}{2}w^{ij}(x)\frac {\partial}{\partial x^i}\wedge \frac
{\partial}{\partial x^j}+ y^aA^{ij}_a(x)\frac {\partial}{\partial
x^i}\wedge \frac {\partial}{\partial y^j}+ \frac {1}{2}
y^ay^bB^{ij}_{ab}(x)\frac {\partial}{\partial y^i}\wedge \frac
{\partial}{\partial y^j} \ .
   \0(2.8)
$$
For a later utilization, we also give

 DEFINITION 2.4. A bivector field $W$ on $TM$ which is locally of
the form $(2.7)$ (respectively, $(2.8)$) is called a {\it
polynomially graded (respectively, graded) bivector field}.

In this case we may speak of a  skew-symmetric    bracket
$$
\{F,G\}_W:=W(dF,dG) \ \ \ (F,G\in C^{\infty}(TM))
  \0(2.9)
$$
which satisfies the Leibniz rule, but, generally, not the Jacobi
identity.

 PROPOSITION 2.5.  {\it If $W$ is a graded Poisson structure on
$TM,$ the equality
$$
\{ l(\alpha ),f \}=-l(D_{df}\alpha) \ , \ \ \alpha \in \Omega
^1(M), \ f\in C^{\infty}(M)
  \0(2.10)
$$
defines a flat contravariant connection on the Poisson manifold
$(M,w).$}

{\it Proof.} By a contravariant connection on $(M,w)$ we
understand a contravariant derivative on the bundle $T^{*}M$ with
respect to the Poisson structure \cite{Va1}.

With $(2.3),$ condition $(2.10)$ means that
$$
D_{df}\alpha :=-{\gamma}_{\alpha}f \ ,
    \0(2.11)
$$
and $(2.4)$ has the equivalent form
$$
D_{df}(h\alpha )=hD_{df}\alpha +((df)^{\sharp}h)\alpha \ , \ \
\alpha \in \Omega ^1(M), \ f,h \in C^{\infty}(M) \ ,
$$
which is the characteristic property of a contravariant connection
 on a Poisson manifold.

 Let us extend $(2.11)$ by
$$
D_{g(df)}\alpha :=gD_{df}\alpha \ , \  \ \ \ \forall g\in
C^{\infty}(M) \ .
   \0(2.12)
$$
The extension is correct because it is compatible with definition
$(2.10)$:  if $g(df)=dh \ $ $(h \in C^{\infty}(M)),$ then $$
-l(D_{dh}\alpha )=\{ l(\alpha),h  \}=W(dl(\alpha ),dh)$$ $$=
gW(dl(\alpha ),df)=g\{ l(\alpha),f\}=-gl(D_{df}\alpha),$$ hence
$D_{dh}\alpha =gD_{df}\alpha$ as needed.

The curvature of this connection is \cite{Va1}
$$
C_{D}(df,dg)\alpha=D_{df}D_{dg}\alpha -D_{dg}D_{df}\alpha
-D_{\{df,dg\}}\alpha  \ ,
$$
and it is easy to see that its annulation is equivalent to the
Jacobi identity
$$
\{\{l(\alpha),f\},g\}+\{\{f,g\},l(\alpha )\}+\{\{g,l(\alpha)\},f\}
\ = \ 0 \ . \ \Endproof
  \0(2.13)
$$

 REMARK 2.6. For any polynomially graded bivector field $W$ such
 that the first term of $(2.7)$ is a Poisson bivector on $M,$ it
 follows similarly that $(2.11)$ and $(2.12)$ define a
 contravariant connection $D$ on $(M,w)$ but, generally,
its curvature  is not  zero.

Now,  let us make some remarks concerning the operator $\Psi $
 of a graded Poisson structure on $TM,$ where
$$
\{l(\alpha),l(\beta )\}= s(\Psi (\alpha ,\beta )) \ .
  \0(2.14)
$$
With $(2.10),$ the second relation   $(2.6)$ becomes
$$
\Psi (\alpha , f\beta )=f\Psi (\alpha , \beta ) -\frac {1}{2}
(D_{df}\alpha \otimes \beta+\beta \otimes D_{df}\alpha  ) \ .
 \0(2.15)
$$
Hence,  $\Psi :T^{*}M\times T^{*}M\longrightarrow {\odot}^2T^{*}M
$ is a bidifferential operator of the first order.

The relation $(2.15)$ allows us to derive
 the local coordinate expression of $\Psi .$
Put
$$
D_{dx^i}dx^j=\Gamma ^{ij}_kdx^k \ , \ \ \ \ \alpha ={\alpha
}_idx^i, \ \beta ={\beta}_jdx^j \ .
  \0(2.16)
$$
It follows that
$$
\Psi (\alpha , \beta )=\alpha _i  \beta _j\Psi (dx^i,dx^j)
  \0(2.17)
$$
$$+\left( {\Gamma}^{kj}_p\delta ^i_q {\beta}_j \frac {\partial
{\alpha }_i}{\partial x^k}- {\Gamma}^{ki}_p\delta ^j_q {\alpha }_i
\frac {\partial {\beta }_j}{\partial x^k} \right)dx^p\odot dx^q
 +w^{kh}\frac {\partial {\alpha }_p}{\partial x^h} \frac
{\partial {\beta }_q}{\partial x^k}dx^p\odot dx^q  \ . $$

 PROPOSITION 2.7. {\it If $G$ is a symmetric covariant tensor
field on $M$ and $\tilde{G}=G_{i_1...i_k}y^{i_1}...y^{i_k}$ is its
corresponding polynomial (see $(1.14))$ then, for any graded
Poisson bivector field $W$ on $TM$, one has  }
$$
\{ \tilde{G},f\}_W=-\widetilde{D_{df}G} \ .
  \0(2.18)
$$

{\it Proof.} Here, $D_{df}$ is the extension of the operator of
contravariant derivative $D$ to $S(T^{*}M)$ i.e.,
$$
(D_{df}G)(X_1,...,X_k)=X_f^{w}(G(X_1,...,X_k))-\sum \limits
_{i=1}^kG(X_1,...,D_{df}X_i,...,X_k) \ ,
$$
where $X_1,...,X_k\in {\mathcal{V}}^1 (M),$ and $D_{df}X$ is
defined by
$$
<D_{df}X, \lambda >=X_f^{w}<X,\lambda>-<X,D_{df}\lambda>, \ \ \
\lambda \in \Omega ^1(M) \ .
$$

 Using the Leibniz rule, we have
$$
\{\tilde{G},f\}=\{G_{i_1...i_k},f\}y^{i_1}...y^{i_k}+ \sum \limits
_{l=1}^k\{y^{i_l},f\}G_{i_1...i_k}y^{i_1}...\hat{y^{i_l}}...y^{i_k}
\ .
$$
But,
$$
\{y^{i},f \}=\{ l(dx^{i}),f  \}=-l(D_{df}dx^{i}   )= -\frac
{\partial f}{\partial x^a}{\Gamma}_h^{ai}y^h \ .
$$
Hence
$$
\{\tilde{G},f\}=-y^{i_1}...y^{i_k}\left( X_f^{w}(G_{i_k...i_k})+
\sum \limits _{l=1}^k  \frac {\partial f}{\partial x^a}{\Gamma
}^{ah}_{i_l}
  G_{i_1...i_{l-1}hi_{l+1}...i_k}  \right) \ .
  $$
The same expression is found for $-\widetilde{(D_{df}G)}.$
 \Endproof

 PROPOSITION 2.8. {\it If we define an operator $D_{df}$
which acts on   $\Psi $ by
$$
(D_{df}\Psi )(\alpha ,\beta ):=D_{df}(\Psi (\alpha , \beta))- \Psi
(D_{df}\alpha , \beta)-\Psi (\alpha ,D_{df} \beta) \ ,
   \0(2.19)
$$
the  Jacobi identity
$$
\{\{l(\alpha),l(\beta)\},f\}+\{\{l(\beta),f\},l(\alpha)\} + \{\{
f,l(\alpha)\},l(\beta)\}=0
  \0(2.20)
$$
is equivalent to}
$$
(D_{df}\Psi )(\alpha ,\beta) =0 \ , \ \ \forall \alpha ,\beta \in
\Omega ^1(M).
  \0(2.21)
$$

{\it Proof. } Express $(2.20)$ by means of $(2.10),$ $(2.14)$ and
$(2.18)$ for $G=\Psi ( \alpha , \beta).$
 \Endproof

Notice that
$$(D_{df}\Psi ) (\alpha ,h \beta)=h(D_{df}\Psi )(\alpha , \beta)-
[C_D(df,dh)\alpha ]\odot \beta \ .
  \0(2.22)
$$
Hence $D_{df}\Psi $ is a bidifferential operator of the second
order. Furthermore, from $(2.22)$ we can see that $(2.21)$ is
invariant by  $\alpha \mapsto f\alpha, \ $  $\beta \mapsto g\beta
$ ($f,g\in C^{\infty}(M)$) iff the curvature $C_D=0.$

In order to discuss the  Jacobi identity
$$
\sum \limits _{(\alpha ,\beta ,\gamma  )} \{ \{l (\alpha),l(\beta)
\}, l(\gamma) \} =0 \
 \0(2.23)
$$
(putting indices between parentheses denotes that summation is on
cyclic permutations of these indices), let us remark the existence
of an operator $\Xi $ such that
$$
\{ s(G), l (\gamma) \}=\widetilde{\Xi (G, \gamma)} \ ,
 \0(2.24)
$$
where $G\in S_2(T^{*}M),$  $\gamma \in \Omega ^1(M),$
 $\Xi (G, \gamma)$ is a symmetric
$3$-covariant tensor field on $M,$ and tilde is the isomorphism
$(1.15)$.

By replacing $G$ by $fG$ and $\gamma $ by $h\gamma ,$ where
$f,h\in C^{\infty}(M),$ we get
$$
\Xi (fG ,h\gamma )=fh\Xi (G ,\gamma )-f(D_{dh}G)\odot \gamma +
hG\odot D_{df}\gamma + \{ f , h \}_{w}G\odot \gamma \ \ ,
 \0(2.25)
$$
which can be used to get the local coordinate expression
$$
\Xi (G, \gamma) = G_{ij}  \gamma _k\Xi (dx^i\odot dx^j,dx^k)+
\frac {1}{3}\sum \limits _{(i,j,k)} ( -G_{hj}\frac {\partial
\gamma _k }{\partial x^a}\Gamma ^{ah}_i
     \0(2.26)
$$
$$ -G_{hi}\frac {\partial \gamma _k }{\partial x^a}\Gamma ^{ah}_j
+\gamma _h\frac {\partial G_{ij} }{\partial x^a}\Gamma ^{ah}_k +
w^{ab} \frac {\partial G_{ij} }{\partial x^a}\frac {\partial
\gamma _k }{\partial x^b} )dx^i\odot dx^j\odot dx^k \ ,
$$
where $\Gamma $ are the local coefficients of the contravariant
connection $D$ defined by $(2.16).$

 Using the operator $\Xi ,$ the Jacobi identity $(2.23)$ becomes
$$
\sum \limits _{(\alpha ,\beta ,\gamma  )}\Xi (\Psi (\alpha
,\beta), \gamma )=0 \ .
 \0(2.27)
$$

We may summarize our analysis by

 PROPOSITION 2.9. {\it The graded bivector field $W$ on $TM$ is
a Poisson bivector iff:
\newline
$a)$ the induced bivector field $w$ on $M$ is Poisson;
\newline
$b)$ the associated contravariant connection $D$ is flat;
\newline
$c)$ the equalities $(2.21)$ and $(2.27)$ hold.

In this case, the projection $\pi :(TM,W)\longrightarrow (M,w)$ is
a Poisson mapping.}

\vspace{0.3cm} To get examples, we consider the following
situation.

Suppose that the symplectic foliation $S$ of an $n-$dimensional
Poisson mani\-fold $(M,w)$ is contained in a regular foliation
${\mathcal{F}}$ on $M$, such that  $T{\mathcal{F}}$ is a {\it
foliated bundle}
 i.e.,  there are local bases
$\{Y_u\}$ $(u=1,...,p ,$ $p= rank \ \mathcal{F})$
 of $T\mathcal{F}$ with  transition  functions constant along
 the leaves of $\mathcal{F}.$
 Consider a decomposition
$$
TM=T{\mathcal{F}}\oplus \nu  \mathcal{F} \ ,
 \0(2.28)
 $$
where $\nu  \mathcal{F}$ is a complementary subbundle of $T
\mathcal{F},$ and  $\mathcal{F}-$adapted local coordinates
$(x^a,y^u) \ \ (a=1,...,n-p) $
 on $M$ \cite{Va2}. Then,
$$\begin{array}{l}
T{\mathcal{F}}=span \left\{\frac {\partial}{\partial
y^u}\right\}=span\{Y_u\},\vspace{2mm}\\
 \nu {\mathcal{F}} =span \left\{X_a:=\frac {\partial}{\partial
 x^a}-t^u_a \frac
{\partial}{\partial y^u}\right\},\end{array}
   \0(2.29)
 $$
 for some local function
$t^u_a=t^u_a(x,y)$. Furthermore, if $\{Y_u \},$ $\{ {\tilde{Y}}_v
\}$ are local bases of the foliated structure of $T{\mathcal{F}}$
over the open neighborhoods $U,$ $ {\tilde{U}} \subseteq M,$ then
$$
{\tilde{Y}}_v=a^u_v(x)Y_u \ , \ \ \ (u,v=1,...,p)
  \0(2.30)
$$
over the connected components of $U\cap\tilde U$.

Since $S\subseteq \mathcal{F},$ the Poisson bivector $w$ is of the
form
 $$
w=\frac {1}{2}w^{uv}(x,y)\frac {\partial}{\partial y^u}\wedge
\frac {\partial}{\partial y^v} \ \ \ \ \ \ (w^{vu}=-w^{uv}).
  \0(2.31)
 $$

Now,  $\forall \  V\in TM$,  $V={\xi}^aX_a+\eta ^uY_u$, and we may
consider   $(x^a,y^u,\xi ^a \ , \eta ^u)$ as {\it distinguished
local  coordinates}  on $TM.$ The transition functions of these
coordinates over the connected components of intersections of
coordinate neighborhoods are of the form
$$
{\tilde{x}}^a ={\tilde{x}}^a(x), \
{\tilde{y}}^u={\tilde{y}}^u(x,y), \ {\tilde{\xi}}^a =
\frac{\partial {\tilde{x}}^a}{\partial x^b}\xi ^b, \
{\tilde{\eta}}^u = b^u_v(x)\eta ^v \ ,
  \0(2.32)
$$
where $b^u_va^v_w= \delta ^u_w$ and $a,b=1,...,n-p;$
$u,v=1,...,p.$

 PROPOSITION 2.10. {\it Under  the previous hypotheses,
the tangent bundle $TM$  has a
 graded Poisson bivector $W,$ which  has the expression $(2.31)$  with respect to the
 distinguished local coordinates.}

 {\it Proof.} It follows from $(2.32)$  that
$$
W=\frac {1}{2}w^{uv}(x,y)\frac {\partial}{\partial y^u}\wedge
\frac {\partial}{\partial y^v}
 \0(2.33)
$$
is a global tensor field on $TM.$ Moreover, since $[W,W]$ has the
same expression as on $M,$ $W$ is  a Poisson bivector.

To prove that $W$ is  graded we also consider   natural
coordinates $({\tilde{x}}^a,{\tilde{y}}^u,z^a,$ $z^u)$ on $TM,$
where $(z^a,z^u)$ are the vector coordinates with respect to the
bases $\{ {\partial}/{\partial {\tilde{x}}^a},
{\partial}/{\partial {\tilde{y}}^u}\}.$ The transition functions
to these coordinates are of the following local form
$$
{\tilde{x}}^a=x^a, \ {\tilde{y}}^u=y^u, \ z^a=\xi ^a, \
z^u=-t^u_a(x,y)\xi ^a+ {\alpha}^u_v(x,y)\eta ^v \ ,
  \0(2.34)
$$
where the coefficients $\alpha ^u_v$ are defined by
$Y_v=\sum_u\alpha ^u_v ({\partial}/{\partial y^u}).$ Accordingly,
$$
\frac {\partial}{\partial y^u}=\frac {\partial}{\partial
{\tilde{y}}^u}+\left(-\frac {\partial t^v_a}{\partial y^u}\xi ^a+
\frac {\partial {\alpha}^v_t}{\partial y^u}\eta ^t\right)\frac
{\partial}{\partial z^v}  \ \ (a=1,...,n-p; \ u,v,t=1,...,p),
$$
and $(2.34)$ shows that $(2.33)$ turns into an expression of type
$(2.8)$. \Endproof

Proposition 2.10 has the following interesting particular cases:

a) The Poisson structure $w$ of $M$ is regular, and the bundle
$TS$ is a foliated bundle; in this case, we take $
{\mathcal{F}}=S.$

b) $S$ is contained in a leaf-wise, locally affine, regular
foliation $\mathcal{F}.$ This means that we have
${\mathcal{F}}-$adapted, local coordinates $(x^a,y^u)$ with local
transition functions
$$
{\tilde{y}}^v=p^v_u(x)y^u+q^v(x) \ ,
$$
and we may use the local vector fields $Y_u= {\partial}/{\partial
y^u} \ .$

c) The Poisson manifold $(M,w)$ has a flat linear connection
$\nabla ,$ possibly with torsion. Then, we may take as leaves of
$\mathcal{F}$ the connected components of $M,$ and the vector
fields $Y_u$ to be local $\nabla -$parallel vector fields. (Then,
in $(2.30)$ we have locally constant coefficients $a^u_v.$)

In particular, the result applies for a locally affine manifold
$M$ (where $\nabla $ has no torsion), and for a parallelizable
manifold $M$ (where we have global vector fields $Y_u$).

As a consequence, we see that Proposition 2.10 holds for the
Lie-Poisson structure of any Lie coalgebra ${\mathcal{G}}^{*}$
\cite{Va1}, which means that $T{\mathcal{G}}^{*}={\mathcal{G}}^{*}
\times {\mathcal{G}}^{*}$ has a graded Poisson structure.
\section*
{\begin{center}{\bf 3. Transversal Poisson structures of
foliations and graded bivector fields on tangent
bundles}\end{center}} The results of the previous section indicate
that the conditions for the existence of a graded Poisson
structure on a tangent bundle $TM$ are rather restrictive. On the
other hand, we will show in this section that more general, but
still interesting, graded bivector fields always exist.

We begin with the following general definition. Let
${\mathcal{F}}$ be an arbitrary regular foliation, with
$p-$dimensional leaves,  on an $n-$dimensional manifold $N.$ We
denote by $C^{\infty}_{fol}(N)$ the  space of differentiable
functions on $N$ which  are constant along the leaves
 of $\mathcal{F}$ ({\it foliated functions}).

DEFINITION 3.1. A {\it transversal Poisson  structure} of
$(N,\mathcal{F})$ is  a bivector field  $w$ on $N$  such that
$$
\{f,g\}:=w(df,dg)  \ , \ f,g\in C^{\infty}(N) \0(3.1)
$$
restricts to a Lie algebra bracket on  $C^{\infty}_{fol}(N).$

 PROPOSITION 3.2. {\it The bivector field  $w\in
{\mathcal{V}}^2(N)$ defines a
   transversal Poisson  structure of the foliation ${\mathcal F}$  iff
 $$
  ({\mathcal{L}}_Yw)|_{Ann
\ T{\mathcal{F}}}=0 \ ,  \ \ \
 [w,w]|_{Ann \ T{\mathcal{F}}}=0 \ ,
   \0(3.2)
$$
for all $Y\in \Gamma (T{\mathcal{F}}).$}

{\it Proof.} The annihilator space $Ann \ T{\mathcal{F}}\subseteq
\Omega ^1(N)$ is
$$
Ann \ T{\mathcal{F}}= \ span\{df  \ / \ f\in C^{\infty}_{fol}(N)
\} \
$$
i.e.,  $f\in C^{\infty}_{fol}(N)$ iff, $\forall \ Y \in
 \Gamma (T{\mathcal{F}}),$ $Yf=0.$

 Accordingly, if  $\ f,g
 \in C^{\infty}_{fol}(N)$ one has
$$
({\mathcal{L}}_Yw)(df,dg)=Y(w(df,dg))=Y\{f,g \} \ ,
$$
and we see that  the first condition $(3.2)$ is equivalent with
  $\{f,g\} \in C^{\infty}_{fol}(N)$, $\forall f,g\in
  C^{\infty}_{fol}(M).$

The second condition $(3.2)$  is a direct consequence of
 the formula (e.g. \cite{YKS1}):
$$
[w,w](df,dg,dh)= 2\sum \limits _{(f,g,h)}\{\{f,g\},h\}  . \
\Endproof
   \0(3.3)
$$

Consider again  a decomposition $(2.28),$ and
$\mathcal{F}-$adapted local coordinates  $(x^a,y^u)$ $ \
(a=1,...,n-p, u=1,...,p)$ on $N$ such that $(2.29)$ holds (with no
reference to any fields $Y_u$ this time). Then
$$
w=\frac {1}{2}w^{ab}X_a\wedge X_b+w^{au}X_a\wedge  \frac
{\partial}{\partial y^u} +\frac 12w^{uv} \frac {\partial}{\partial
y^u} \wedge  \frac {\partial}{\partial y^v} \ ,
   \0(3.4)
$$
and the first condition  $(3.2)$ means that, locally,
$w^{ab}=w^{ab}(x).$

Although this is not our main subject, we will derive some more
facts about transversal Poisson structures of foliations.

 PROPOSITION 3.3. {\it The Hamiltonian vector field $X_f:=i(df)w$
of a foliated function $f$ is a foliated vector field (i.e.,
projectable  on the space of leaves).}

{\it Proof.}  A vector field $Z\in \Gamma TN$ is foliated if
$\forall \ Y\in \Gamma (T\mathcal{F}),$ ${\mathcal{L}}_YZ\in
\Gamma (T\mathcal{F}).$  But, if $Y\in \Gamma (T\mathcal{F})$ and
$g\in C^{\infty}_{fol}(N)$ then, by $(3.2),$
 $$
 ({\mathcal{L}}_YX_f)g=dg({\mathcal{L}}_Y({\sharp}_w(df) )=
 ({\mathcal{L}}_Yw)(df,dg)+w(d(Yf),dg)=0 \ .
 $$
Therefore ${\mathcal{L}}_YX_f \in \Gamma (T\mathcal{F}).$
\Endproof

 DEFINITION 3.4. The generalized distribution $\mathcal{D}$ defined by
$$
{\mathcal{D}}_x=span \{ \ Y(x),X_f(x)  \ /  \ Y\in \Gamma
(T{\mathcal{F}}), \ f\in C^{\infty}_{fol}(N)\} \ \ (x\in N)
   \0(3.5)
$$
is called the {\it characteristic distribution of $w$ on
$(N,{\mathcal{F}}).$}

 PROPOSITION 3.5. {\it The characteristic distribution
${\mathcal{D}}$ of a transversal Poisson structure of a foliation
is completely integrable, and each leaf $\Sigma$ of $\mathcal{D}$
is a presymplectic manifold,  with a presymplectic $2-$form of
kernel  $ {T\mathcal{F}}|_{\Sigma}.$ }

{\it Proof.} Brackets of the form $[Y_1,Y_2],$ $[Y,X_f],$
$Y_1,Y_2, Y \in \Gamma (T\mathcal{F}),$ $ f\in
C^{\infty}_{fol}(N)$ belong to $\mathcal{D}$ because $\mathcal{F}$
is a  foliation, and because of Proposition 3.3. The latter also
shows that $\forall \ f,g,h \in C^{\infty}_{fol}(N),$
$$
dh([X_f,X_g]-X_{ \{f,g\} })=0 \ ,
$$
whence
$$
[X_f,X_g]=X_{ \{f,g\} }+Y, \ \ Y\in \Gamma (T\mathcal{F}) \ .
$$
Thus, the distribution $\mathcal{D}$ is involutive.

Furthermore, let $U$ be an $\mathcal{F}-$adapted coordinate
neighborhood, and $p:U\longrightarrow V$, $V:=U/{U\cap
\mathcal{F}}$ the submersion onto the corresponding space of
slices. Because of Proposition 3.3, the distribution
$p_{*}({\mathcal{D}})$ exists on $V,$ and, obviously, it precisely
is tangent to the symplectic distribution of the Poisson structure
induced by the first term of $(3.4)$ on $V.$ It follows that
$p_{*}(\mathcal{D})$ has a constant dimension along the integral
paths of the vector fields $p_{*}X_f$ $(f\in
C^{\infty}_{fol}(N)).$ Hence
${\mathcal{D}}=p^{-1}_{*}(p_{*}({\mathcal{D}}))$ has a constant
dimension along the integral paths of the vector fields $X_f.$
$\mathcal{D}$ also has a constant dimension along the integral
paths of vector fields $Y\in \Gamma (T\mathcal{F})$ because
 $p_{*}(\mathcal{D})$ does not change along such paths.

 Now, the complete integrability of $\mathcal{D}$ follows from one of the
 versions of the Frobenius-Sussmann-Stefan theorem, Theorem
 2.9$''$
 of \cite{Va1}.

 The leaves $\Sigma$ of the characteristic distribution $\mathcal{D}$
 are immersed submanifolds of $N$ which are foliated by the
 corresponding restriction of $\mathcal{F},$ and are sent by the
 submersion $p:U\longrightarrow V:=U/_{{U\cap \mathcal{F}}}$
 encountered above to symplectic manifolds, included in the
 symplectic leaves, say $\sigma ,$ of the projection of the first
 term of $(3.4).$ It is obvious that the symplectic forms of $\sigma$
 lift to a global, closed $2-$form $\lambda $ on $\Sigma ,$
with the kernel $ {T\mathcal{F}}|_{\Sigma}.$ \Endproof

As a matter of fact, we may notice that $w$ produces more than
just a presymplectic structure on the leaves $\Sigma $ of
$\mathcal{D}.$  It also produces the generalized distribution
$$
E:={\sharp}_w \ Ann(T{\mathcal{F}})=span \ \{X_f \ / \ f\in
C^{\infty}_{fol}(N)\}
$$
which has a restriction of constant rank on each leaf $\Sigma ,$
 such that $T\Sigma =T({\mathcal{F}}|_{\Sigma})\oplus E|_{\Sigma}
 \ .$

\vspace{3mm} Now we return to the tangent bundles $TM.$ All of
them have the vertical foliation $\mathcal{F}$  by fibers with the
tangent distribution $V:=T\mathcal{F}.$

The set of  foliated functions on $TM,$ may be identified with
$C^{\infty}(M).$

 PROPOSITION 3.6. {\it  Any polynomially  graded bivector field
 $W$ on $TM,$ which is $\pi $ related with a Poisson structure of $M$ is
 a transversal Poisson structure of $(TM,V).$    }

{\it Proof.} $\pi $ is the projection $TM\longrightarrow M,$ and
if we take $W$ as in $(2.7),$ $W$ is  $\pi -$related with the
tensor $w$  defined on $M$ by the first term of $(2.7).$ Then,
$(3.1)$ obviously holds. \Endproof

 DEFINITION 3.7. A transversal Poisson  structure of
 the vertical foliation  of $TM$ will be called a
{\it semi-Poisson structure} on $TM.$

In particular, the structures of Proposition 3.6 are polynomially
graded semi-Poisson structures.

\vspace{0.3cm}

In what follows, we will discuss a class of graded semi-Poisson
 structures of a tangent bundle $TM$
 and show how to construct all the  graded semi-Poisson bivector fields
 on $TM$ which have a given induced  Poisson
 structure $w$ on the base manifold.

 Let us consider a Poisson bivector $w$ on $M.$ Recall that a
semispray (a second order differential equation) \cite{Leon} on
$M$ is a vector field $S$ on $TM$ such that $FS=E,$ where $F=
({\partial }/{\partial y^i})\otimes dx^i$ is the natural almost
tangent structure and $E=y^i( {\partial }/{\partial y^i}) $ is the
Euler vector field on $TM.$  The local coordinate expression of
$S$ is of the form
$$
S=y^i\frac {\partial }{\partial x^i} +{\sigma}^i(x,y)\frac
{\partial }{\partial y^i} \ .
 \0(3.6)
$$

Let $\nabla $ be a torsionless linear connection on $M,$ with the
local coefficients ${\Gamma}^k_{ij}$ and $S$ its associated
semispray (the  geodesic spray) given by
$$
S=y^i\frac {\delta }{\delta x^i} \ , \ \ \frac {\delta }{\delta
x^i}:= \frac {\partial }{\partial x^i}-y^k\Gamma^j_{ik}\frac
{\partial }{\partial y^j} \ .
$$

 PROPOSITION 3.8. {\it If $(M,w)$ is a Poisson manifold then
 the bivector field
$$
W=-\frac 12{\mathcal{L}}_Sw^C \ ,
 \0(3.7)
$$
where $w^C$ is the complete lift of $w$ to $TM,$  defines a
graded semi-Poisson structure on $TM.$ }

{\it Proof.} If the local coordinate expression of $w$ is $(1.1),$
$w^C$ is given by $(1.2),$ and we get
$$
W=\frac {1}{2}w^{ij}(x)\frac {\partial}{\partial x^i}\wedge \frac
{\partial}{\partial x^j}-y^aw^{ik}{\Gamma}^j_{ka}\frac
{\partial}{\partial x^i}\wedge \frac {\partial}{\partial
y^j}-\frac 14 y^ay^b (\frac {{\partial}^2w^{ij}}{\partial x^a
\partial x^b}
\0(3.8)
$$
 $$
-  \frac {\partial w^{ij}}{\partial x^k}{\Gamma}^k_{ab}+
w^{kj}\frac {\partial {\Gamma}^i_{ab}}{\partial x^k}- w^{ki}\frac
{\partial {\Gamma}^j_{ab}}{\partial x^k} +2\frac {\partial
w^{kj}}{\partial x^b}{\Gamma}^i_{ka}-2 \frac {\partial
w^{ki}}{\partial x^b}{\Gamma}^j_{ka} ) \frac {\partial}{\partial
y^i}\wedge \frac {\partial}{\partial y^j} \ . \  \Endproof
$$

From  $(3.7),$ it follows that
$$
\{F_1,F_2\}_W:=W(dF_1,dF_2)=
 \0(3.9)
$$
$$
=-\frac {1}{2}({\mathcal{L}}_S\{F_1,F_2\}_{w^C}
-\{{\mathcal{L}}_SF_1,F_2\}_{w^C}-\{F_1,{\mathcal{L}}_SF_2\}_{w^C})
\ , \
$$
where $F_1,F_2 \in C^{\infty}(TM).$

For further reference, we will say that $W$ of $(3.7),$ $(3.8)$ is
the {\it graded $\nabla -$lift } of the Poisson structure $w$ of
$M$. We are going to describe it in a different form below.

First, $\forall H\in S_k(T^{*}M),$ define $^s {\nabla } H\in
S_{k+1}(T^{*}M)$ by
$$
^s {\nabla } H(X_1,...,X_{k+1})=\frac {1}{k+1}\sum \limits
_{i=1}^{k+1} (\nabla _{X_i}H)(X_1,...,\hat{X_i}...X_{k+1}) \ .
$$
Then, with the notation of $(1.15),$ it follows easily that
 $$
{\mathcal{L}}_S\tilde{H}=\widetilde{^{s} {\nabla } H} \ .
 \0(3.10)
 $$

If $G_1,G_2\in S(T^{*}M),$ using  $(1.16)$ and $(3.10),$  we get
the explicit formula
$$
\{{\tilde{G}}_1, {\tilde{G}}_2   \}_W=-\frac 12\iota(^{s} {\nabla
}<G_1,G_2>
    \0(3.11)
$$
$$
-<^{s} {\nabla } G_1,G_2>-<G_1,^{s} {\nabla } G_2>) \  ,
$$
where $< \ , \ >$ is the Schouten-Nijenhuis bracket of symmetric
tensors, and $\iota$ is the isomorphism $(1.15).$

 PROPOSITION 3.9. {\it The graded $\nabla -$lift  $W$ of $w$
 is characterized by:

$i)$ The Poisson structure induced by $W$ on the base manifold
coincides with the given Poisson structure $w$ on $M,$ i.e.
$$
\{f,g\}_W=\{f,g\}_{w} \ , \  \ \forall f,g \in C^{\infty}(M);
 \0(3.12)
$$

$ii)$ for every $f \in C^{\infty}(M)$ and $\alpha \in \Omega
^1(M)$
$$
\{l(\alpha),f\}_W=-l(\nabla _{X_f}\alpha ) \ ;
 \0(3.13)
$$

$iii)$ for any Pfaff forms $\alpha $ and $\beta $ of $M$ we have
$$
\{l(\alpha),l(\beta)\}_W=-\frac 12 \iota(^{s} {\nabla } < \alpha ,
\beta>-<^{s} {\nabla } \alpha ,\beta >-<\alpha ,^{s} {\nabla }
\beta >) \ .
 \0(3.14)
$$
}

{\it Proof.} $i)$ If $f,g \in C^{\infty}(M),$ from $(3.9)$ and the
definition of $w^C$ in Section 1, we get
$$
\{f,g\}_W= \frac 12 (\{{\mathcal{L}}_Sf,g \}_{w^C}+ \{{f,
\mathcal{L}}_Sg\}_{w^C})=\frac 12
(\{l(df),g\}_{w^C}-\{l(dg),f\}_{w^C})
$$
$$
=\frac 12(X^w_fg-X^w_gf)=\{ f,g \}_w \ .
$$

$ii)$ For $f\in C^{\infty}(M)$ and  $\alpha \in \Omega ^1(M)$,
 $(3.11)$ becomes
$$
\{l(\alpha ) ,f\}_W= -\frac {1}{2}({\mathcal{L}}_S\{l(\alpha
),f\}_{w^C}-\{{\mathcal{L}}_S l(\alpha ),f\}_{w^C}-\{l(\alpha
),{\mathcal{L}}_Sf \}_{w^C}) \ .
 \0(3.15)
$$
Here, we have
$$\{l(\alpha ),f\}_{w^C}=-\alpha (X_f) \ , \ \
 {\mathcal{L}}_S\{l(\alpha ),f\}_{w^C}=-l(d(\alpha (X_f))) \ ,
$$
and, with
 $(3.10)$, $(1.13)$ and $(1.16)$,
$$
\{{\mathcal{L}}_S l(\alpha ),f\}_{w^C}= -l(i_{X_f}( ^s {\nabla }
\alpha )) \ .
$$
Finally, we have
$$
\{l(\alpha ),{\mathcal{L}}_Sf \}_{w^C}=\{l (\alpha  ), l(df)
\}_{w^C}=l(\{\alpha ,df   \})
$$
$$
=-l({\mathcal{L}}_{X_f}\alpha )=-l(d(\alpha (X_f))+i_{X_f}d\alpha
) \ .
$$

With these results $(3.15)$ gives
$$
\{l (\alpha ),f\}_W= -\frac 12 l[i_{X_f}(d \alpha +  ^s {\nabla }
\alpha )] \ .
  \0(3.16)
$$
Since  the torsion of $\nabla $ vanishes
 we have
 $$
2(d\alpha )(X,Y)=(\nabla _X\alpha )Y-(\nabla _Y\alpha )X \ , \ \
X,Y \in \chi (M) \ ,
 $$
and
$$
d \alpha +  ^s {\nabla } \alpha =\nabla \alpha \ ,
$$
where $\nabla \alpha $ is the $2-$covariant tensor field defined
by $\nabla \alpha (X,Y)=(\nabla _X\alpha )(Y), \ $ $ \forall
X,Y\in {\mathcal{V}}^1(M), $ and $(3.16)$  exactly becomes
$(3.13).$

$iii)$  $(3.14)$ is a direct consequence of $(3.11).$  \Endproof

 REMARK 3.10. Comparing the relation $(3.13)$ with $(2.10)$ we see
that the contravariant derivative $D$ associated to the graded
semi-Poisson structure $W$ is  the contravariant derivative
induced by the  linear connection $\nabla $ (see \cite{Va1}).

 REMARK 3.11. The relation $(3.14)$ provides us the expression of the
operator $\Psi _W$ associated to $W$ (see $2.14):$
$$
{\Psi}_W(\alpha ,\beta )= -\frac 12(^s{\nabla} <\alpha ,\beta >-<
^s{\nabla } \alpha ,\beta >- <\alpha ,  ^s{\nabla }\beta >) \ .
 \0(3.17)
$$

 \vspace{0.3cm}

Now, we will prove

 PROPOSITION 3.12. {\it Let $(M,w)$ be a Poisson manifold. The
graded semi-Poisson structures $W$ on $TM$ for which the canonical
projection $\pi :(TM,W)\longrightarrow (M,w)$ is a Poisson mapping
are defined by the relations
$$
\{f,g\}_W=\{f,g\}_{w} \ , \ \{l(\alpha),f\}_W=-l(D_{df}\alpha ),
\;
 \{l(\alpha),l(\beta)\}_W =s(\Psi (\alpha ,\beta )) \
,
$$
$f,g \in C^{\infty}(M),$  $\alpha , \beta \in \Omega ^1(M),$ where
$D$ is an arbitrary contravariant connection of $(M,w)$ and the
operator $\Psi $ is given by
$$
\Psi ={\Psi}_0 +A+T \ ,
 \0(3.18)
$$
with terms as follows: ${\Psi}_0$ is the
  operator $\Psi $ of a fixed  graded semi-Poisson structure
  $W_0,$
$A:T^{*}M\times T^{*}M\longrightarrow {\odot}^2T^{*}M$ is a
skew-symmetric, first order, bidifferential operator with the
 property
$$
A(\alpha ,f\beta )=f A(\alpha , \beta)-{\tau}(df , \alpha )\odot
\beta \ ,
 \0(3.19)
$$
where $\tau $ is a  tensor field of type $(2,1)$ on $M, $
 and $T\in \Gamma(({\wedge}^2 TM)\otimes ({\odot}^2T^{*}M)).$}

{\it Proof.} If $D$ is the contravariant derivative associated to
$W$ in Remark 2.6, then, to change it,  means to use  a connection
$D^{'}=D+\tau ,$ where $\tau $ is a  tensor field of type $(2,1)$
on $M.$ Accordingly, from $(2.15)$ it follows that $
{\Psi}^{'}-\Psi $
 is a bidifferential operator with the property $(3.19).$
Then, with the contravariant connection $D$ chosen, we see from
$(2.15)$ again, that the only possible change of $\Psi $ consist
in adding a tensor $T.$ \Endproof

 REMARK 3.13. An example of operator ${\Psi}_0$ is
 provided by
${\Psi}_W$ given by $(3.17).$

Notice that a given Poisson structure $w$ on $M$ may have no
graded Poisson lift to $TM.$ In particular, a flat contravariant
connection $D$ may not exist.  Indeed \cite{Va1}, one can mimic
 the  Chern-Weil construction of characteristic classes
  and  associate to each Poisson manifold $(M,w),$
Pontriagin-Poisson classes $p_k(M,w)$ which are the image of the
usual Pontriagin classes in the $LP-$cohomology
 by the homomorphism
 $$
 \sharp : [\lambda ]\in H^k_{DR}(M)\longmapsto [{\lambda}^{\sharp}]\in
H^k_{LP}(M) \ .
$$
 If a flat $D$ exists, all $p_k(M,w)=0.$ Thus, if a non zero
Pontriagin-Poisson class exists, there is no flat
  connection $D.$
\section*{\begin{center}{\bf 4. Horizontal  lifts of a Poisson
structure}
\end{center}}
In this section, we define {\it horizontal lifts} of a Poisson
bivector $w$ to the tangent  bundle of the Poisson manifold
$(M,w)$
 and study the
conditions for these lifts to be Poisson bivectors, and to be
compatible with the complete lift $w^C.$

 Let $M$ be an $n-$dimensional
 manifold and  $\pi :TM\longrightarrow M$ its tangent bundle.
On $TM,$ we consider a {\it nonlinear connection}  i.e., a
distribution $\mathrm{H}$, called {\it horizontal,}  such that
$T(TM)=\mathrm{H}\oplus \mathrm{V}, $ where $\mathrm{V}$ denotes
the vertical distribution tangent to the fibers of $TM$
\cite{{Leon},{Vop}}. If $(x^i)$  are local coordinates on $M$ and
 $(x^i,y^j)$ ($(i,j=1,...,n)$) are the induced coordinates on $TM$ (see Section 1),
 we have bases of the form
$$
\mathrm{V}=span \left\{  \frac {\partial}{\partial y^i}  \right\}
\ , \ \mathrm{H}=span \left\{  \frac {\delta}{\delta x^i}:=\frac
{\partial}{\partial x^i}-\Gamma ^j_i\frac {\partial}{\partial y^j}
\right\} \ ,
 \0(4.1)
$$
and ${\Gamma}^i_j$ are called the  {\it coefficients of the
 connection.}

Equivalently, the nonlinear connection may be seen as an almost
product structure $\Gamma $ on $TM$ such that
 the eigendistribution   corresponding to the eigenvalue $-1$ is the vertical
 distribution ${\mathrm{V}}$  \cite{Leon}.
Then
 $$
h=\frac 12(Id +\Gamma ):TM\longrightarrow \mathrm{H}
 $$
is the {\it horizontal projector}  of $\Gamma ,$ and the {\it
curvature} $R $ of the connection  is the {\it Nijenhuis tensor}
$$
R(X,Y)=-N_h(X,Y)=-[hX,hY]+h[hX,Y]+h[X,hY]-h[X,Y] \ ,
$$
where $X,Y \in {\mathcal{V}}^1 (TM).$ $R$ vanishes if at least one
argument is in $V$, and always takes values in $V$, hence,
locally, we may write \cite{Leon}
$$
R=\frac 12R^k_{ij}dx^i\wedge dx^j\otimes \frac {\partial}{\partial
y^k} \ ,\hspace{5mm} R^k_{ij}=\frac {\delta {\Gamma}^k_j}{\delta
x^i}-\frac {\delta {\Gamma}^k_i}{\delta x^j} \ .\0(4.2)
$$

Then, we get $$\left[\frac {\delta}{\delta x^i},\frac
{\delta}{\delta x^j}\right]=-R^k_{ij} \frac {\partial}{\partial
y^k} \ ,\hspace{5mm}\left[\frac {\delta}{\delta x^i},\frac
{\partial}{\partial y^j} \right]=\frac {\partial
{\Gamma}^k_i}{\partial y^j}
 \frac {\partial}{\partial y^k} \ .
 \0(4.3)
$$
In particular, $H$ is involutive iff $R=0$.

Let us consider a bivector $w$ on the base manifold $M,$ having
the local coordinate expression $(1.1).$

 DEFINITION 4.1. The {\it horizontal lift} of $w$  to the tangent
bundle $TM$, with respect to the connection $\Gamma$ is  the
(global) bivector field $w^H$ defined by
$$
w^H=\frac 12 w^{ij}(x)\frac {\delta}{\delta x^i}\wedge \frac
{\delta}{\delta x^j} \ .
 \0(4.4)
$$

Notice that the horizontal lift $(4.4)$ is different from that of
\cite{KY}.

  PROPOSITION 4.2. {\it Let $(M,w)$ be a Poisson manifold. If the
  horizontal distribution $H$
  is defined by a linear connection $\nabla $ on
$M$,  the bivector $w^H$  defines a graded semi-Poisson structure
on $TM.$}

{\it Proof.} With respect to the bases
 $({\partial}/{\partial x^i},
 {\partial}/{\partial y^j})$, the expression
  of $w^H$ is of the form $(2.8).$
\Endproof

PROPOSITION 4.3. {\it A horizontal lift  $w^H$ is a Poisson
bivector on  $TM$ iff $w$ is a Poisson bivector on the base
manifold $M$ and}
$$
R(X_f^H,X_g^H)=0, \ \  \ \forall f,g \in C^{\infty }(M) \ ,
 \0(4.5)
$$
where $X_f$ denotes the $w-$Hamiltonian vector field of $f$ and
$X_f^H$ is the horizontal lift of  $X_f$  \cite{KY}.

{\it Proof.}
 A straightforward computation yields
 the  Schouten-Nijenhuis bracket
$$
[w^H,w^H]=\frac 13 \left(\sum \limits _{(i,j,k)}w^{hk} \frac
{\partial w^{ij}}{\partial x^h} \right)\frac {\delta}{\delta x^i}
\wedge \frac {\delta}{\delta x^j}\wedge \frac {\delta}{\delta
x^k}+
 \0(4.6)
$$
$$
+w^{hi}w^{lj}R^k_{hl}\frac {\delta}{\delta x^i}\wedge \frac
{\delta}{\delta x^j}\wedge \frac {\partial}{\partial y^k} \ .
$$
Since the vanishing of the first term of $(4.6)$
 is equivalent to $[w,w]=0$ on $M,$
 $w^H$ is a Poisson bivector on $TM$ iff $w$ is a Poisson
bivector on $M$ and
$$
w^{ih}w^{jl}R^k_{hl}=0 \ .
 \0(4.7)
$$
The latter equation  has the  equivalent form
$$
R({({\sharp}_w \alpha)}^H,{({\sharp}_w \beta)}^H)=0, \ \ \forall
\alpha, \beta \in {\Omega}^1(M) \ ,
 \0(4.8)
$$
which is also  equivalent to $(4.5).$ \Endproof

 REMARKS 4.4.

$i)$ If $w^H$ is a Poisson bivector, the projection $\pi :(TM,w^H)
\longrightarrow (M,w)$ is a Poisson mapping.

$ii)$  If $w$ is defined by a symplectic form on $M,$  condition
$(4.5)$ becomes $R=0.$

 COROLLARY 4.5. {\it If $(M,w)$ is a Poisson manifold and the
connection $\Gamma $ on $TM$ is defined by a covariant derivative
$\nabla $ on $M$,   the bivector $w^H$  defines a
  Poisson structure on $TM$ iff the curvature $C_D$ of the contravariant
  connection induced by $\nabla$ on $TM$  vanishes.
  In this case, $w^H$ is a graded Poisson structure on $TM$.}

{\it Proof.} Remember that $D$ is defined by
$D_{df}={\nabla}_{X_f},$ and we may see this operator as acting
either on $T^{*}M$ or on $TM$ \cite{Va1}.

If ${\Gamma }^k_{ij}$ are the connection coefficients of $\nabla
,$  ${\Gamma }^k_i=\Gamma ^k_{ij}(x)y^j$ and $
R^k_{ij}=y^hR^k_{hij} \ ,$ where $R^k_{hij}$ are the components of
the curvature $R_{\nabla}.$ Condition $(4.7)$ becomes
$$
R_{\nabla }(\sharp \alpha , \sharp \beta)Z=0, \ \ \forall \alpha ,
\beta \in {\Omega}^1(M), \ \forall Z\in {\mathcal{V}}^1 (M) \ ;
 \0(4.9)
$$
(equivalently,
$$
R_{\nabla }(X_f,X_g)Z =0  \ \ \forall f,g \in C^{\infty}(M) \ , \
\forall Z\in {\mathcal{V}}^1 (M) \ ).
 \0(4.9')
$$
This condition is equivalent to $C_D=0.$ \Endproof.

If the connection $\Gamma $ on $TM$ is defined by a covariant
derivative $\nabla $ on $M,$ the conditions for the graded
bivector field  $w^H$ to be Poisson are simpler than those of
Proposition 2.9, and $(2.21)$ , $(2.27)$ must be consequences of
the conditions of Proposition 4.3. Furthermore, one can check that
 the operators $\Psi $ and $\Xi$ of $w^H$ (see (2.14), (2.24))
are given by
$$
{\Psi }_{w^H}(\alpha ,\beta)(X,Y)=\frac 12[w(({\nabla}_{.}\alpha
)X, ({\nabla}_{.}\beta )Y)+ w(({\nabla}_{.}\alpha )Y,
{\nabla}_{.}\beta )X)] \ ,
   \0(4.10)
$$
$$
\Xi (G,\gamma)(X,Y,Z)=\frac {1}{3!}\sum \limits
_{(X,Y,Z)}w({\nabla}_{.}G)(X,Y),({\nabla}_{.}\gamma)Z) \ ,
   \0(4.11)
$$
where $\alpha,\beta,\gamma\in\Omega^1(M)$, $G\in S_2(T^*M)$, and
$\nabla.$ means that we create a $1$-form which is evaluated on
$Z\in{\mathcal V}^1(M)$ by the application of $\nabla_Z$.

\vspace{2mm} Let us consider  an arbitrary Poisson structure $W$
on $TM.$

Following \cite{Va3}, we would like to know whether there are
semisprays on $TM$ which are Hamiltonian vector fields with
respect to $w^H.$

 PROPOSITION 4.6. {\it If the Poisson bivector $w$ on $M$ is not
defined by a symplectic  structure, there are no $w^H-$Hamiltonian
semisprays on $TM.$}

{\it Proof.} If $F\in C^{\infty}(TM),$ then
$$
X_F^{w^H}=w^{ij}\frac {\delta F}{\delta x^i}\frac
{\partial}{\partial x^j}- w^{ik} \frac {\delta F}{\delta
x^i}\Gamma ^j_k \frac {\partial}{\partial y^j}
$$
and $(3.6)$ shows that  $X_F^{w^H}$ is a semispray iff
$$
w^{ij}\frac {\delta F}{\delta x^i} =y^j \ .
 \0(4.12)
$$
$(4.12)$ implies
$$
-w^{jh}\frac {\partial}{\partial y^k}\left(\frac {\delta F}{\delta
x^h} \right)=\delta ^j_k \ ,
$$
therefore, $(w^{jh})$ is a non singular matrix.  \Endproof

\vspace{3mm}Recall that two Poisson structures on a manifold $M$
are {\it compatible} if the bivector fields $w_1$ and $w_2$
satisfy the condition
$$
[w_1,w_2]=0 \ ,
 \0(4.13)
$$
or,  equivalently,   $w_1+w_2$  also is a Poisson bivector field.

If $w^H$ is a Poisson bivector, it is natural to discuss its
compatibility with the complete lift $w^C$ of $w.$

  PROPOSITION 4.7. {\it Let $w$ be a Poisson structure,  and $\nabla$ a
 symmetric linear connection on $M$ such that the associated contravariant
 connection of $TM$ has zero curvature. Then the Poisson bivector
 $w^H$ is compatible with the complete lift $w^C$ iff
 $$
i_{X_f}(\nabla ^2w)=0, \ \ \forall f \in C^{\infty}(M) \ ,
 \0(4.14)
 $$
where $\nabla ^2 w=\nabla \nabla w$ is the tensor field of type
$(2,2)$ on $M$ defined by }
$$
(\nabla ^2w)(X,Y)=(\nabla _X(\nabla w))Y=\nabla _X\nabla _Yw-
{\nabla }_{\nabla _XY}w, \ \ X,Y \in {\mathcal V}^1 (M) \ .
$$

{\it Proof. } We  compute the bracket $[w^H,w^C]$  using the
auxiliary notations
$$
a^{ij}:=y^k \frac {\partial w^{ij}}{\partial x^k}+\Gamma
^i_kw^{kj}- \Gamma ^j_kw^{ki} \ , \ \ t^k_{hl}=\frac {\partial
\Gamma ^k_h}{\partial y^l}- \frac {\partial \Gamma ^k_l}{\partial
y^h} \ .
 \0(4.15)
$$

By a straightforward computation, it follows that $[w^H,w^C]=0$ is
equivalent to
$$
w^{hi}w^{jl}t^k_{hl}=0 \ , \ \ w^{hi}\left( \frac {\delta
a^{jk}}{\delta x^h}-a^{lj} \frac {\partial \Gamma ^k_h}{\partial
y^l}+a^{lk} \frac {\partial \Gamma ^j_h}{\partial y^l}\right)=0 \
.
 \0(4.16)
$$

If  $\Gamma $ comes from a symmetric linear connection $\nabla $
on $M,$ the first condition $(4.16)$ holds, and the second
condition (4.16) is the coordinate expression of (4.14). \Endproof

 REMARK 4.8. For any $w\in {\mathcal V}^2(M)$, one can see that
 $Q=w^H+w^C$ is a Poisson bivector iff
$w$ and $w^H$ are Poisson bivectors and $w^C$ is compatible with
$w^H.$ \vspace{3mm}\\

{\bf Acknowledgement} During the work on this paper, Gabriel
Mitric held a postdoctoral grant at the University of Haifa,
Israel. He wants to thank the University of Haifa, its Department
of Mathematics, and, personally, Professor Izu Vaisman for the
grant and for hospitality during the period when the postdoctoral
program was completed.

\begin{center}
{\small \begin{tabular}{ll}
Gabriel Mitric&Izu Vaisman\\
Catedra de Geometrie&Department of Mathematics\\
Universitatea ``Al. I. Cuza" Ia\c{s}i\hspace*{2.5cm}&University of
Haifa\\Rom\^{a}nia&Israel\\
gmitric@uaic.ro\hspace{7cm}&vaisman@math.haifa.ac.il
\end{tabular}}
\end{center}
\end{document}